\documentclass[11pt]{article}
\usepackage{color,latexsym,amsthm,amsmath,amssymb,path,url,enumitem}

\newcommand{\ignore}[1]{}

\newcommand{\parag}[1]{\vspace{2mm}

\noindent{\bf #1} }

\newtheorem{theorem}{Theorem}[section]
\newtheorem{lemma}[theorem]{Lemma}

\newcounter{problems}
\newtheorem{problem}[problems]{Problem}

\newsavebox{\smallProofsym}                     
\savebox{\smallProofsym}                        %
        {
        \begin{picture}(7,7)                    %
        \put(0,0){\framebox(6,6){}}             %
        \put(0,2){\framebox(4,4){}}             %
        \end{picture}                           %
        }                                       %
\usepackage{graphicx,psfrag}
\usepackage[rflt]{floatflt}

\newcommand{\placefig}[2]
        {\includegraphics[width=#2]{#1.eps}}
\usepackage{dsfont}

\newcommand{\ZZ}{\ensuremath{\mathbb Z}}
\newcommand{\RR}{\ensuremath{\mathbb R}}

\newcommand{\pts}{\mathcal P}

\def\hd{{\hat{D}}}

\newcommand{\s}{{\mathsf{subset}}}
\DeclareMathOperator*{\EE}{\mathbb{E}}
\newcommand{\lattice}{\mathcal L}

\def\eps{{\varepsilon}}

\begin{document}
\pagenumbering{arabic}

\title{Distinct Distances: Open Problems and Current Bounds}

\author{
Adam Sheffer\thanks{Department of Mathematics, Baruch College, City University of New York, NY, USA.
{\sl adamsh@gmail.com}. Supported by NSF grant DMS-1710305} }

\date{}

\maketitle

\begin{abstract}
We survey the variants of Erd\H os' distinct distances problem and the current best bounds for each of those.
\end{abstract}

\section{Introduction}
Given a set $\pts$ of $n$ points in $\RR^2$, let $D(\pts)$ denote the number of distinct distances that are determined by pairs of points from $\pts$. Let $D(n) = \min_{|\pts|=n}D(\pts)$;
that is, $D(n)$ is the minimum number of distinct distances that a set of $n$ points in $\RR^2$ can determine.
In his celebrated 1946 paper \cite{erd46}, Erd\H os derived the bound $D(n) = O(n/\sqrt{\log n})$, and conjectured that this bound is tight.
More specifically, Erd\H os showed that a $\sqrt{n} \times \sqrt{n}$ section of the integer lattice $\ZZ^2$ determines $\Theta(n/\sqrt{\log n})$ distinct distances.
Though over 70 years have passed since Erd\H os considered this lattice structure, no configuration that determines an asymptotically smaller number of distinct distances was discovered.

For the celebrations of his 80th birthday, Erd\H os compiled a survey of his favorite contributions to mathematics \cite{erd96}, in which he wrote
 \begin{quotation}\em ``My most striking contribution to geometry is, no doubt, my problem on the number of distinct distances. This can be found in many of my papers on combinatorial and geometric problems."
 \end{quotation}

After over sixty years and a series of increasingly larger lower bounds, Guth and Katz \cite{GK15} derived the bound $D(n) = \Omega(n/\log n)$, almost matching the current best upper bound. 
A comprehensive study of the previous bounds can be found in \cite{GIS11}.\footnote{See also William Gasarch's webpage: \url{http://www.cs.umd.edu/~gasarch/erdos_dist/erdos_dist.html}} 
To derive their bound, Guth and Katz developed several novel techniques, relying on tools from algebraic geometry, 19th century analytic geometry, and more.
Notice that a small gap of $O(\sqrt{\log n})$ remains between the current best lower and upper bounds.
\begin{problem} \label{pr:main}
Find the exact asymptotic value of $D(n)$.
\end{problem}

Since Problem \ref{pr:main} is almost completely solved, one might wonder what is the purpose of this survey.
This problem is just one out of many challenging distinct distances problems, most of which are still wide open (and also originally posed by  Erd\H os).
For some of these problems, such as the ones presented in Section \ref{sec:structure}, hardly anything non-trivial is known after decades of work.
The study of distinct distances problems is an active sub-field, with a constant stream of new results.
This survey is an attempt to keep track of this progress.

For readers who are only interested in the main open problems, it is the personal view of the author (and likely of others) that currently the most challenging/interesting distinct distances problems are:
\begin{itemize}
\item Finding the minimum number of distinct distances spanned by $n$ points in $\RR^d$. See Problem \ref{pr:Rd}.
\item Characterizing the point sets in $\RR^2$ that span a small number of distinct distances. See Section \ref{sec:structure}.
\end{itemize}

The survey is partitioned into sections according to sub-families of distinct distances problems. 
Section \ref{sec:structure} discusses the structure of planar point sets that span few distinct distances. 
Section \ref{sec:rest} surveys problems in $\RR^2$ in which the point set is restricted in some manner. 
Section \ref{sec:HighD} considers distinct distances in $\RR^d$. 
Section \ref{sec:Bipartite} studies bipartite problems.
Section \ref{sec:subset} discusses subsets of point sets where no distance repeats more than once.
Section \ref{sec:kl} is about using local distance properties to derive global distance properties.
Section \ref{sec:additive} studies problems that are related to Additive Combinatorics.
Finally, Section \ref{sec:additional} contains a few problems that do not fit into any of the other sections.

\parag{Acknowledgements.} The author is indebted to the people who helped improving this survey: Adrian Dumitrescu, William Gasarch, Ben Lund, Cosmin Pohoata, Micha Sharir, and Frank de Zeeuw.

\section{The structure of point sets with few distinct distances} \label{sec:structure}

In this section we discuss the characterization of point sets in $\RR^2$ that span few distinct distances.
After decades of studying this topic, hardly anything is known about it.
One might say that this family of problems is the one for which we know the least, and we may still not have the correct tools for handling it.

Since the asymptotic value of $D(n)$ is still unknown, we consider sets $\pts$ of $n$ points in $\RR^2$ that satisfy $D(\pts)=O(n/\sqrt{\log{n}})$, and refer to such sets as \emph{near-optimal}.
All the point sets in this section are planar.

\begin{problem} \label{pr:Structure}
Characterize the near-optimal point sets.
\end{problem}

\begin{figure}[h]

\centerline{\placefig{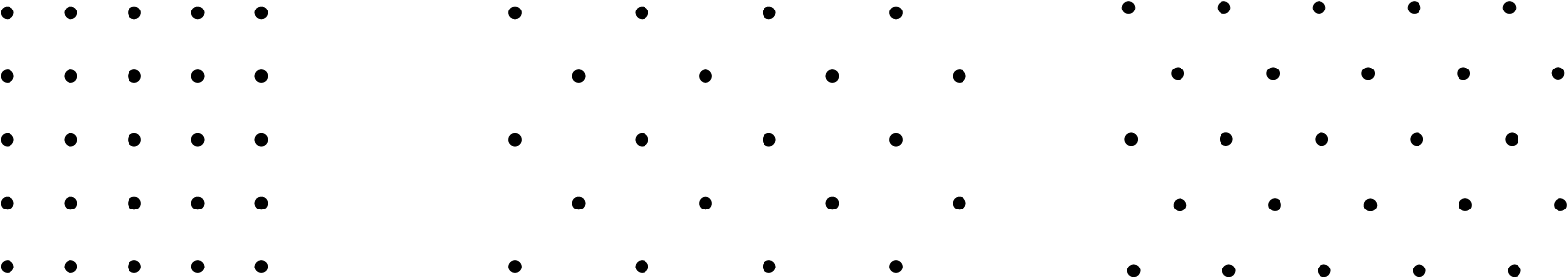}{0.9\textwidth}}
\vspace{-1mm}

\caption{\small \sf (a) An section of $\ZZ^2$.
(b) A lattice that can be obtained from (a) either by rotation and uniform scaling, or by removing every point whose coordinates sum to an odd number.
(c) A triangular lattice.}
\label{fi:grids}
\vspace{-2mm}
\end{figure}

Erd\H os asked whether every near-optimal set ``has lattice structure" \cite{erd86}.
To make this question more rigorous, we first consider some of the known near-optimal sets.
In the introduction we already mentioned that a $\sqrt{n} \times \sqrt{n}$ section of the integer lattice $\ZZ^2$ determines $O(n/\sqrt{\log n})$ distinct distances (e.g., see Figure \ref{fi:grids}(a)). This was observed by Erd\H os, who noticed that this is an immediate corollary of the following theorem from number theory.
\begin{theorem} \label{th:LandRam} {\bf (Landau-Ramanujan \cite{BR95,BMO11,Land08})}
The number of positive integers smaller than $n$ that are the sum of two squares is $\Theta(n/\sqrt{\log n})$.
\end{theorem}
Every distance in the $\sqrt{n} \times \sqrt{n}$ integer lattice is the square root of a sum of two squares between 0 and $n$. Thus, Theorem \ref{th:LandRam} implies that the number of distinct distances in this case is $\Theta(n/\sqrt{\log n})$.

The above implies that the $\sqrt{n} \times \sqrt{n}$ integer lattice is a near-optimal set.
More generally, for any integer $c \ge 1$, every $n$ point subset of the $c\sqrt{n} \times c\sqrt{n}$ integer lattice is a near-optimal set, since we can still apply Theorem \ref{th:LandRam} in such cases. 
We may obtain additional near-optimal sets by applying operations such as translations, rotations, and uniform scalings.
For example, the lattice in Figure \ref{fi:grids}(b) can be obtained either by rotating and scaling the $\sqrt{n} \times \sqrt{n}$ section of $\ZZ^2$, or by removing from the $\sqrt{2n} \times \sqrt{2n}$ section of $\ZZ^2$ every point whose coordinates sum to an odd number.

We can further generalize the above. For any integer $r>1$, we define the rectangular lattice
\[ {\mathcal L}_r = \{ (i,j\sqrt{r}) \ \mid \ i,j\in\ZZ \quad \text{and} \quad  1 \le i,j \le \sqrt{n} \}. \]
Every set ${\mathcal L}_r$ spans $\Theta(n/\sqrt{\log n})$ distinct distances (e.g., see \cite{Sblog14}).
By applying various transformations to these rectangular lattices, we get additional near-optimal lattices. 
For example, the triangular lattice, which corresponds to the vertices in a tiling of equilateral triangles (e.g., see Figure \ref{fi:grids}(c)), can be obtained by taking a rectangular lattice and removing every other vertex.
Erd\H os and Fishburn \cite{EF96} conjectured that, for infinitely many values of $n$, there is an $n$-point subset of the triangular lattice that minimizes the number of distinct distances (not only asymptotically)

Hardly anything is known regarding Erd\H os's conjecture that every near-optimal set has a lattice structure.
The author of this survey suggests that perhaps the near-optimal sets are exactly the ones that could be obtained from the sets ${\mathcal L}_r$.
As a first step, Erd\H os \cite{erd86} suggested to determine whether every near-optimal point set contains $\Omega(\sqrt{n})$ points on a line, and thus most of the set can be covered by a small number of lines. Since this also appears to be quite difficult, Erd\H os asked whether there exists a line with $\Omega(n^\varepsilon)$ points of the set. Embarrassingly, even this weaker variant remains open.

\begin{problem} \label{pr:OnLine} {\bf (Erd\H os \cite{erd86})}
Prove or disprove: For a sufficiently small $\varepsilon>0$, every near-optimal point set contains $\Omega(n^\varepsilon)$ points on a common line.
\end{problem}

It is known that for every near-optimal set $\pts$ of $n$ points, there exists a line $\ell$ such that $|\ell\cap\pts| = \Omega(\log n)$ (e.g., see \cite{Sblog14b}).
Moreover, it is shown in \cite{LSdZ15} that for every near-optimal set $\pts$ of $n$ points and $0 < \alpha \le 1/4$, either there exists a line or a circle that contains $n^{\alpha}$ points of $\pts$, or there exist $n^{8/5-12\alpha/5-\eps}$ distinct lines that contain $\Omega(\sqrt{\log n})$ points of $\pts$.
Thus, one possible approach for solving Problem \ref{pr:OnLine} might be to prove that for any near-optimal set $\pts$ with  many lines that contain $\Omega(\sqrt{\log n})$ points of $\pts$, there exists a line containing $\Omega(n^{\eps})$ points of $\pts$.

Sheffer, Zahl, and de Zeeuw \cite{SZZ13} considered the complement problem --- proving that no line can contain many points of a near-optimal set. 
They proved
that for every near-optimal set $\pts$ of $n$ points, every line contains  $O(n^{7/8})$ points of $\pts$. This bound was recently improved to  $O(n^{43/52})$ in \cite{RRNS15}.
In \cite{SZZ13}, it is also proved that for every near-optimal set $\pts$ of $n$ points, every \emph{circle} contains  $O(n^{5/6})$ points of $\pts$.
Pach and de Zeeuw \cite{PZ13} showed that if a set of $n$ points is contained in a constant-degree curve $\gamma$, then these points span $\Omega(n^{4/3})$ distinct distances, unless $\gamma$ contains a line or a circle.
That is, for every near-optimal set $\pts$ of $n$ points, any constant degree algebraic curve that does not contain lines and circles contains  $O(n^{3/4})$ points of $\pts$.
By combining these three results, we obtain that for every near-optimal set $\pts$ of $n$ points, every constant-degree algebraic curve contains  $O(n^{43/52})$ points of $\pts$.

\begin{problem} \label{pr:OptimalCurves}
Prove or disprove: For every near-optimal set $\pts$ of $n$ points and $\varepsilon>0$, every constant-degree curve contains $O(n^{0.5+\varepsilon})$ points of $\pts$.
\end{problem}

More problems related to the structure of point sets with few distances can be found in Section \ref{sec:additive}.

\section{Restricted point sets in $\RR^2$} \label{sec:rest}

\begin{figure}[h]

\centerline{\placefig{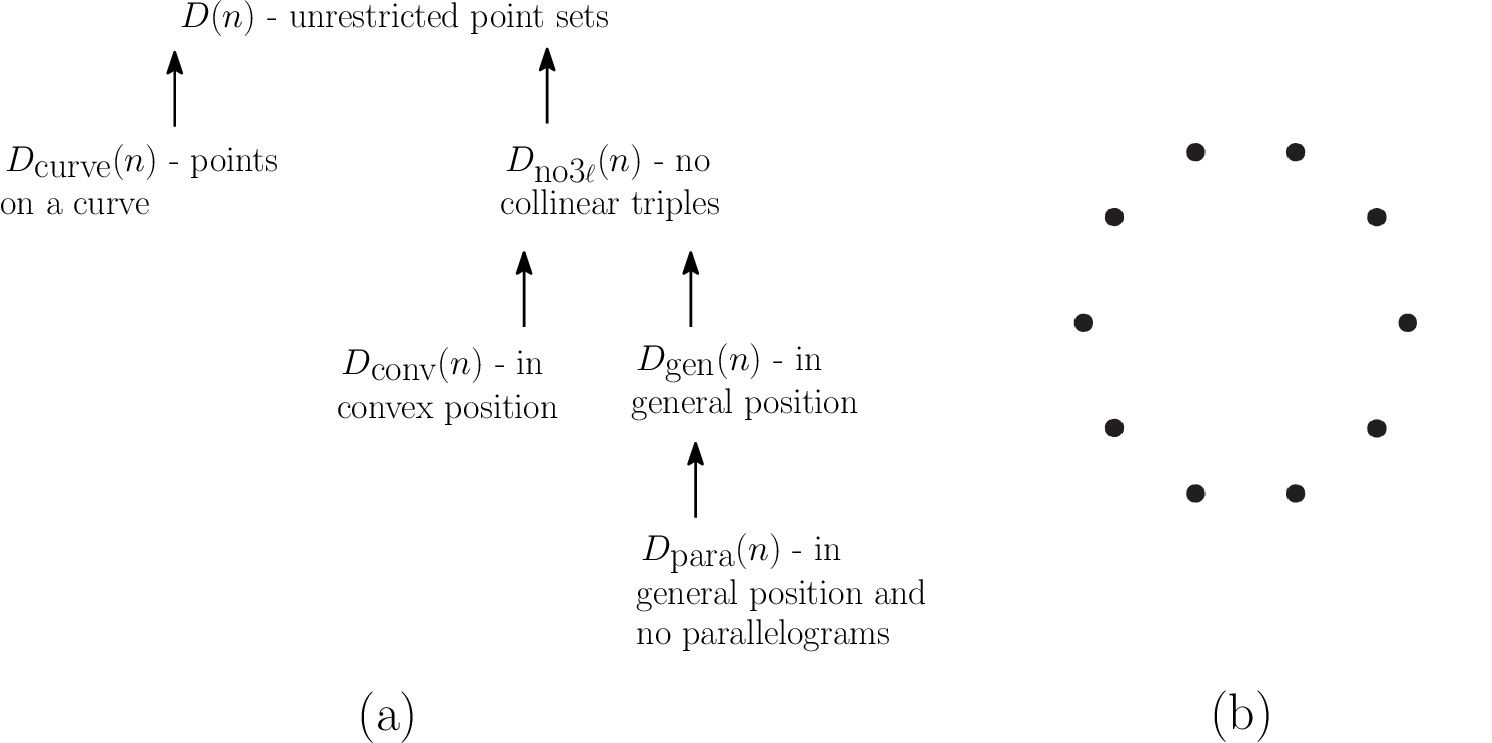}{0.92\textwidth}}
\vspace{-1mm}

\caption{\small \sf (a) The hierarchy of the restricted point sets in $\RR^2$. 
Every arrow goes from a problem to a less restricted generalization. (b) The vertices of a regular $n$-gon lie on a common circle.}
\label{fi:sec2}
\vspace{-2mm}
\end{figure}

In this section we consider variants of the planar distinct distances problem where the point sets are restricted in some way. 
The current best bounds for these problems are listed in Table \ref{ta:restR2}; see Figure \ref{fi:sec2} and the text below for an explanation of the notation used in the table.
All of the point sets in this section are in $\RR^2$.

\begin{table}[h]
\begin{center}
\caption{The current best bounds for restricted sets of points in $\RR^2$.\label{ta:restR2}}
\vspace{2mm}

\begin{tabular}{|c|c|c|}
\hline
{\textbf {Variant}} & {\textbf {Lower bound}} & {\textbf {Upper bound}} \\
\hline\hline
$D_{\text{curve}}(n)$  & $\Omega(n^{4/3})$ \cite{PZ13} & $\displaystyle O(n^2)$ (trivial) \\ \hline
$D_{\text{no3}\ell}(n)$  & $\displaystyle \lceil (n-1)/3\rceil$ (Szemer\'edi) & $\displaystyle \lfloor n/2\rfloor$ \cite{erd46} \\ \hline
$\hd_{\text{conv}}(n)$  & $\left(\frac{13}{36}+\frac{1}{22701}\right)n +O(1)$ \cite{Dum06,NPPZ11} & $\displaystyle \lfloor n/2\rfloor$ \cite{erd46} \\ \hline
$D_{\text{gen}}(n)$  & $\Omega(n)$ (trivial) & $\displaystyle n2^{O(\sqrt{\log n})}$ \cite{EFPR93} \\ \hline
$\hd_{\text{gen}}(n)$  & $\lceil (n-1)/3 \rceil$ (trivial) & $n-O(1)$ (trivial) \\ \hline
$D_{\text{para}}(n)$  & $\Omega(n)$ (trivial) & $\displaystyle O(n^2/\sqrt{\log n})$ \cite{Dum08} \\ \hline

\end{tabular}
\end{center}
\vspace{-\baselineskip}
\end{table}

For our first problem, we restrict the point set $\pts$ to be contained in some constant-degree algebraic curve $\gamma$.
That is, $\gamma$ is the set of points on which a constant-degree polynomial vanishes.
When $\gamma$ is a line, we can obtain $D(\pts) = n-1$ by taking the points of $\pts$ to be evenly spaced on $\gamma$.
When $\gamma$ is a circle, we can have $D(\pts) = \lfloor n/2\rfloor$ by taking the points of $\pts$ to be the vertices of a regular $n$-gon (see Figure \ref{fi:sec2}(b)).
We denote by $D_{\text{curve}}(n)$ the minimum number of distinct distances spanned by $n$ points on a constant-degree curve $\gamma$, when $\gamma$ does not contain any lines or circles.
Pach and de Zeeuw \cite{PZ13} proved that $D_{\text{curve}}(n) = \Omega(n^{4/3})$.
It seems plausible that the actual value of $D_{\text{curve}}(n)$ is very close to $n^2$.

\begin{problem} \label{pr:curve}
Find the asymptotic value of $D_{\text{\em curve}}(n)$.
\end{problem}

We move to study a family of problems that ask for exact bounds, rather than asymptotic ones. 
That is, problems where the goal is to find the best constant of proportionality. 
Denote by $D_{\text{no3}\ell}(n)$ the minimum number of distinct distances determined by a set of $n$ points, no three of which are collinear.
That is, $D_{\text{no3}\ell}(n) = \min_{|\pts|=n}D(\pts)$, where the minimum is taken over all sets of $n$ points containing no three collinear points.
Notice that the vertices of a regular $n$-gon, such as the set depicted in Figure \ref{fi:sec2}(b), satisfy this property and determine $\left\lfloor \frac{n}{2}\right\rfloor$ distinct distances.
We thus have that $D_{\text{no3}\ell}(n) \le \left\lfloor \frac{n}{2}\right\rfloor$. 
The current best lower bound, due to Szemer\'edi (communicated by Erd\H os in \cite{Erd75}), is $D_{\text{no3}\ell}(n) \ge \left\lceil \frac{n-1}{3}\right\rceil$.
Szemer\'edi also conjectured that $D_{\text{no3}\ell}(n) = \left\lfloor \frac{n}{2}\right\rfloor$; see \cite{erd87,EHP89}.
\begin{problem} \label{pr:No3}
Find the exact value of $D_{ \text{\em no3}\ell}(n)$.
\end{problem}

Szemer\'edi's proof is so simple and elegant that it is hard to resist stating it here.

\begin{lemma}\label{le:no3}
$D_{\text{\em no3}\ell}(n) \ge \lceil \frac{n-1}{3}\rceil$
\end{lemma}
\begin{proof}
Consider a set $\pts$ of $n$ points, no three of which are collinear. 
Let $x$ denote the minimum number satisfying that every point $p\in \pts$ determines at most $x$ distinct distances with the points of $\pts\setminus\{p\}$.

Let $T = \{(a,p,q)\in \pts^3 \mid |ap|=|aq|\}$, where $a,p,q$ are three distinct points and where $(a,p,q)$ and $(a,q,p)$ are counted as the same triple. The proof is based on double counting $|T|$, and we begin by deriving an upper bound for it. Given a pair of points $p,q\in\pts$, the triplet $(a,p,q)$ is in $T$ if and only if $a$ is on the perpendicular bisector of the segment $pq$. By the assumption, each such perpendicular bisector contains at most two points of $\pts$, which implies
\begin{equation}\label{eq:SzUp}
|T| \le 2 \binom{n}{2} = n(n-1).
\end{equation}
For the lower bound, notice that for every point $p\in\pts$, the points of $\pts\setminus\{p\}$ are contained in at most $x$ concentric circles around $p$. We denote these circles as $C_{p,1},\ldots,C_{p,x}$ and set $n_{p,i} = |C_{p,i} \cap \pts|$. Notice that $\sum_{i=1}^xn_{p,i} = n-1$ for every $p\in\pts$. By the Cauchy-Schwarz inequality, we have $\sum_{i=1}^x n_{p,i}^2 \ge \frac{1}{x}(n-1)^2$. This in turn implies
\begin{align} |T| = \sum_{p\in\pts}\sum_{i=1}^x \binom{n_{p,i}}{2} = \frac{1}{2}\sum_{p\in\pts}\sum_{i=1}^x (n_{p,i}^2-n_{p,i}) &\ge \frac{1}{2}\sum_{p\in\pts} \left( \frac{1}{x}(n-1)^2 - (n-1)\right) \nonumber \\
&= \frac{n(n-1)(n-1-x)}{2x}. \label{eq:SzLw}
\end{align}
Combining \eqref{eq:SzUp} and \eqref{eq:SzLw} immediately implies the assertion of the lemma.
\end{proof}

Although the last progress made for Problem \ref{pr:No3} was several decades ago, more recent advances have been obtained for more restricted variants.
Let $D_{\text{conv}}(n)$ denote the minimum number of distinct distances determined by a set of $n$ points in (strict) convex position. Let $\hd_{\text{conv}}(n)$ denote the maximum number satisfying that for any set $\pts$ of $n$ points in convex position, there exists a point $p\in \pts$ such that there are at least $\hd_{\text{conv}}(n)$ distinct distances between $\pts$ and $\pts\setminus\{p\}$.

By considering the regular $n$-gon once again, we get that $D_{\text{conv}}(n) \le \left\lfloor \frac{n}{2}\right\rfloor$ and that $\hd_{\text{conv}}(n) \le \left\lfloor \frac{n}{2}\right\rfloor$. 
Already in his 1946 paper, Erd\H os \cite{erd46} conjectured that $D_{\text{conv}}(n) = \left\lfloor \frac{n}{2}\right\rfloor$. 
This was proven by Altman \cite{Alt63,Alt72}, which led Erd\H os to suggest the stronger conjecture $\hd_{\text{conv}}(n) = \left\lfloor \frac{n}{2}\right\rfloor$.
It is not difficult to verify that Lemma \ref{le:no3} proves the existence of $\left\lceil \frac{n-1}{3}\right\rceil$ distinct distances from a single point.
Since there are no three collinear points in a set in convex position, this lemma also implies that $\hd_{\text{conv}}(n) \ge \left\lceil \frac{n-1}{3}\right\rceil$.
In 2006, Dumitrescu \cite{Dum06} derived the improved bound $\hd_{\text{conv}}(n) \ge \lceil \frac{13n-6}{36}\rceil$. Recently, the slightly improved bound $\hd_{\text{conv}}(n) \ge \left(\frac{13}{36}+\frac{1}{22701}\right)n +O(1)$ was obtained by Nivasch, Pach, Pinchasi, and Zerbib \cite{NPPZ11}.

\begin{problem} \label{pr:conv}
Find the exact value of $\hd_{\text{\em conv}}(n)$.
\end{problem}

We say that a set of points is in \emph{general position} if no three points are collinear and no four points are cocircular. Denote by $D_{\text{gen}}(n)$ the minimum number of distinct distances determined by a set of $n$ points in general position. The convex $n$-gon configuration is not in general position, and it is not known whether $D_{\text{gen}}(n) = \Theta(n)$ or not. The current best upper bound $D_{\text{gen}}(n) = n2^{O(\sqrt{\log n})}$ was derived by Erd\H os, F\"uredi, Pach, and Ruzsa \cite{EFPR93}. This bound is obtained by  considering a very different construction: taking an integer grid $G$ in a $d$-dimensional space (where $d$ is roughly $\sqrt{\log n}$), considering a subset $G'$ of the points of $G$ that lie on a common hypersphere, and projecting $G'$ on a generic plane. 
The hypersphere and the generic projection guarantee that the resulting set  is in general position, while the integer grid structure implies a relatively small number of distinct distances.
For an easy lower bound, note that $D_{\text{gen}}(n) \ge D_{\text{no3}}(n) = \Omega(n)$.

\begin{problem} \label{pr:gen}
Find the asymptotic value of $D_{\text{\em gen}}(n)$.
\end{problem}

\ignore{ -------------------------------------------------ignore-----------------------------------------------
Erd\H os \cite{erd87} also suggested to study the maximum number $\hd_{\text{\em gen}}(n)$ satisfying that for any set $\pts$ of $n$ points in general position, there exists a point $p\in \pts$ such that $\{p\}\times\pts$ determines at least $\hd_{\text{gen}}(n)$ distinct distances. Consider a point set $\pts$ in general position and a point $p\in \pts$. If $p$ has a distance of $d$ from four points of $\pts\setminus\{p\}$, then these four points are on the circle centered at $p$ and of radius $d$, contradicting the general position assumption. Thus, a trivial lower bound is $\hd_{\text{\em gen}}(n) \ge \lceil (n-1)/3 \rceil$. 
No non-trivial lower or upper bounds are known for $\hd_{\text{\em gen}}(n)$, and when discussing the problem, Erd\H os \cite{erd87} wrote ``It is rather frustrating that I got nowhere with\ldots".

\begin{problem} \label{pr:gen}
Find the exact value of $\hd_{\text{\em gen}}(n)$.
\end{problem}

} 

The point configuration that implies $D_{\text{gen}}(n) = n2^{O(\sqrt{\log n})}$ spans many duplicate vectors. This led to denoting by $D_{\text{para}}(n)$ the minimum number of distinct distances determined by a set of $n$ points in general position that do not determine any parallelograms.
Erd\H os, Hickerson, and Pach \cite{EHP89} asked whether $D_{\text{para}}(n)= o(n^2)$. 
This was confirmed by Dumitrescu \cite{Dum08}, who proved $D_{\text{para}}(n)= O(n^2/\sqrt{\log n})$.
For prime $n$, Dumitrescu considered the point set
\[ \left\{ (i,j) \mid i=0,1,\ldots,(n-1)/4, \ j= i^2 \hspace{-2mm}\mod n \right\}. \]
No lower bound better than the trivial $D_{\text{para}}(n)= \Omega(n)$ is known.

\begin{problem} \label{pr:para}
Find the asymptotic value of $D_{\text{\em para}}(n)$.
\end{problem}

\section{Higher dimensions} \label{sec:HighD}

In this section we consider higher-dimensional variants of the distinct distances problem.
Denote by $D_d(n)$ the minimum number of distinct distances that a set of $n$ points in $\RR^d$ can determine. 
As in the planar case, the current best lower bound is obtained by considering an even section of the integer lattice.
That is, in the $d$-dimensional case we consider an $n^{1/d} \times n^{1/d} \times \cdots \times n^{1/d}$ section of $\ZZ^d$.
Every distance in this configuration is the square root of a sum of $d$ squares, each with a value between 0 and $n^{2/d}$. 
Since every positive integer can be written as a sum of four squares, and a large portion of the integers can be written as a sum of three squares, the number of distinct distances that are determined by such a lattice is $O(n^{2/d})$.
That is, for $d\ge 3$ we have $D_d(n) = O(n^{2/d})$. This bound was already observed by Erd\H os in his 1946 paper \cite{erd46}, and is conjectured to be tight.

Solymosi and Vu \cite{SV08} derived the following recursive relations on $D_d(n)$.
\begin{theorem} \label{th:SV} {\bf (Solymosi and Vu \cite{SV08})}

\noindent (i) If $D_{d_0}(n) = \Omega(n^{\alpha_0})$, then for all $d>d_0$, we have
\[D_d(n) = \Omega\left( n^{\frac{2d}{(d+d_0+1)(d-d_0)+2d_0/\alpha_0}} \right). \]
(ii) If $D_{d_0}(n) = \Omega(n^{\alpha_0})$, then for all $d>d_0$ where $d-d_0$ is even, we have
\[D_d(n) = \Omega\left( n^{\frac{2(d+1)}{(d+d_0+2)(d-d_0)+2(d_0+1)/\alpha_0}} \right). \]
\end{theorem}
Recall that $D_2(n) = \Omega(n/\log{n})$. Combining this with Theorem \ref{th:SV}(i) implies\footnote{In the $\Omega^*(\cdot)$ notation we neglect polylogarithmic factors; although the logarithm in the denominator of the lower bound for $D_{2}(n)$ does not exactly fit the formulation of Theorem \ref{th:SV}, the proof remains valid.} $D_3(n) = \Omega^*(n^{3/5})$, while the above lattice example implies $D_3(n) = O(n^{2/3})$. These are the current best bounds for $D_3(n)$.
The current best bounds for larger values of $d$ are obtained by combining Theorem \ref{th:SV}(ii) with the bounds $D_2(n)=\Omega^*(n)$ and $D_3(n)=\Omega^*(n^{3/5})$ as base cases. That is, for even $d\ge4$ we have $D_d(n) = \Omega^*\left(n^{\frac{2d+2}{d^2+2d-2}}\right)$ and for odd $d\ge5$ we have $D_d(n) = \Omega^*\left(n^{\frac{2d+2}{d^2+2d-5/3}}\right)$.
Note that as $d$ goes to infinity $D_d(n)$ approaches the conjectured bound $\Theta(n^{2/d})$.

\begin{problem} \label{pr:Rd}
Find the asymptotic value of $D_d(n)$.
\end{problem}

It seems possible that the techniques that were used by Guth and Katz \cite{GK15} for analyzing $D(n)$ could also be applied to the higher dimensional variant. 
Recently, Bardwell-Evans and Sheffer \cite{BES17} reduced the distinct distances problem in $\RR^d$ into an incidence problem with well-behaved $(d-1)$-flats in $\RR^{2d-1}$.
Deriving the conjectured bound for this incidence problem would settle Problem \ref{pr:Rd}.

Let $D_d^{o}(n)$ denote the minimum number of distinct distances that a set of $n$ points on a hypersphere in $\RR^d$ can determine. Tao \cite{Tao11} observed that the bound from \cite{GK15} remains valid when the point set is on a sphere in $\RR^3$ (or on a hyperbolic plane). That is, $D_3^{o}(n) = \Omega(n/\log{n})$. For a lower bound, place a set of $n$ points on a circle that is on the sphere, so that they form the vertices of a regular planar $n$-gon (recall figure \ref{fi:sec2}(b)). This implies $D_3^{o}(n) = O(n)$.

\begin{problem} \label{pr:sphere}
Find the asymptotic value of $D_3^o(n)$.
\end{problem}

Erd\H os, Fu\"redi, Pach, and Ruzsa \cite{EFPR93} proved $D_4^o(n) = O(n/\log \log n)$ and $D_d^o(n) = O(n^{2/(d-2)})$ for $d>4$.
These bounds are obtained by taking an integer lattice in $\RR^d$ and then choosing a hypersphere that contains many lattice points. 
No lower bound is known beyond the trivial $D_d^o(n) \ge D_d(n)$.

\begin{problem} \label{pr:sphereD}
Find the asymptotic value of $D_d^o(n)$ for $d\ge4$.
\end{problem}

\paragraph{Restricted point sets.}
Charalambides \cite{Chara13} considered the case where a set of $n$ points is contained in a constant-degree curve $\gamma$ in $\RR^d$.
Note that this is the $d$-dimensional variant of Problem \ref{pr:curve}.
Charalambides showed that when $\gamma$ contains an \emph{algebraic helix}, the points may determine only $O(n)$ distinct distances (see \cite{Chara13} for a description of algebraic helices. In $\RR^2$ and $\RR^3$ the only algebraic helices are lines and circles).
On the other hand, if $\gamma$ does not contain any algebraic helices,
the points on it determine $\Omega(n^{5/4})$ distinct distances. 
Raz \cite{Raz16} improved this bound to $\Omega(n^{4/3})$ distinct distances, matching the current best planar bound.
We denote by $D_{\text{\em curve}}^{(d)}(n)$ the minimum number of distinct distances that are determined by $n$ points on constant-degree curve in $\RR^d$ that does not contain an algebraic helix.

\begin{problem} \label{pr:curveD}
Find the asymptotic value of $D_{\text{\em curve}}^{(d)}(n)$.
\end{problem}

It is natural to ask what happens when the points are restricted to a surface or to an algebraic variety of any dimension. 
Sharir and Solomon \cite{ShSo17} studied the case of $n$ points on a constant-degree algebraic surface in $\RR^3$ that contains no planes and no spheres.
In this case, they proved that the number of distinct distances is $\Omega(n^{7/9-\eps})$.
The known techniques yield bounds in several other restricted cases in $\RR^d$.
However, it is not yet clear what the main problems and difficulties are.
For example, in the problem studied by Sharir and Solomon, is the restriction about planes and spheres necessary?
We thus end this section with a deliberately vague problem.

\begin{problem} 
Derive non-trivial distinct distances bounds for points on varieties in $\RR^d$, when $d\ge 3$.
\end{problem}

\section{Bipartite problems} \label{sec:Bipartite}

\begin{figure}[h]

\centerline{\placefig{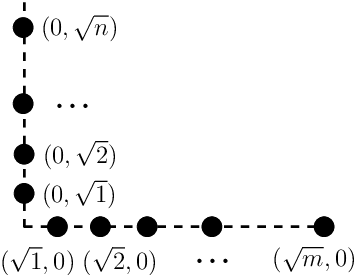}{0.3\textwidth}}
\vspace{-1mm}

\caption{\small \sf While there are few distinct distances between points on different lines, there are many distinct distances between points on the same line.}
\label{fi:OrthLines}
\vspace{-2mm}
\end{figure}

In a \emph{bipartite} distinct distances problem we have two sets of points $\pts_1$ and $\pts_2$, and consider only distances between pairs of points in $\pts_1\times \pts_2$.
That is, we do not care about distances between pairs of points from the same set.
We denote this number of distinct distances as $D(\pts_1,\pts_2)$.
The values of $D(\pts_1)$ and $D(\pts_2)$ may be significantly larger than $D(\pts_1,\pts_2)$.
For example, let $\pts_1$ be a set of $m$ points on the $x$-axis and let $\pts_2$ be a set of $n$ points on the $y$-axis, as depicted in Figure \ref{fi:OrthLines}.
Then $D(\pts_1,\pts_2) = \Theta(m+n)$, $D(\pts_1) = \Theta(m^2)$, and $D(\pts_2) = \Theta(n^2)$.

Let $D(m,n)$ denote minimum number of distinct distances that are determined by point sets in $\RR^2$ of respective sizes $m$ and $n$.
That is, $D(m,n) = \min_{|\pts_1|=m \atop |\pts_2|=n}D(\pts_1,\pts_2)$.
Without loss of generality, we assume that $m\le n$.
We have the trivial upper bound $D(m,n) \le D(m+n) = O(n/\sqrt{\log n})$.
Moreover, Elekes \cite{Elek95} proved that $D(m,n) = O(m^{1/2}n^{1/2})$ when $n\ge 4m^3$.
On the other hand, the lower bound of Guth and Katz bound does not immediately extend to the bipartite case.

\begin{problem} \label{pr:bipartite}
Find the asymptotic value of $D(m,n)$.
\end{problem}

In some sense, Problem \ref{pr:bipartite} asks to extend the Guth and Katz analysis to the bipartite case.
One might expect such an extension to lead to the bound $D(m,n) = \Omega\left(m^{1/2}n^{1/2}/\sqrt{\log n}\right)$.

We next consider bipartite problems where the point sets are restricted to curves.
Let $\ell_1$ and $\ell_2$ be two lines in $\RR^2$.
Let $\pts_1$ be a set of $m$ points on $\ell_1$ and let $\pts_2$ be a set of $n$ points on $\ell_2$.
As illustrated in Figure \ref{fi:OrthLines}, when the two lines are orthogonal we may have $D(\pts_1,\pts_2) = \Theta(m+n)$.
It is not difficult to verify that this bound still holds when the two lines are parallel. 
Purdy conjectured that when the lines are neither parallel nor orthogonal, the number of distinct distances should be superlinear  (e.g., see \cite[Section 5.5]{BMP05}).
We denote as $D_{\text{lines}}(m,n)$ the minimum number of distinct distances in such a scenario.

Elekes and R\'onyai \cite{ER00} proved Purdy's conjecture, though without deriving any specific superlinear lower bound. 
The current best bound, derived in \cite{SS13}, is $D_{\text{lines}}(m,n) = \Omega(\min\{n^{2/3}m^{2/3},m^2,n^2\})$.
The first term in the minimum is the interesting one --- the other two dominate only when one point set is significantly larger than the other.
Elekes \cite{Elek99} observed the upper bound $D_{\text{lines}}(n,n) = O(n^2/\sqrt{\log n})$.

\begin{problem} \label{pr:lines}
Find the asymptotic value of $D_{\text{\em lines}}(m,n)$.
\end{problem}

Among other reasons, Problem \ref{pr:lines} is considered interesting since it has many generalizations, including to problems that do not involve distances (for example, see \cite{ER00,RSS14}).
Improving the known bounds for Problem \ref{pr:lines} tends to lead to improvements for the various generalizations.
Quoting Hilbert \cite{Reid70}: ``The art of doing mathematics is finding that special case that contains all the germs of generality.''

One can generalize Problem \ref{pr:lines} by replacing the lines $\ell_1$ and $\ell_2$ with constant-degree algebraic curves.
Specifically, let $\pts_1$ be a set of $m$ points on a curve $\gamma_1$ and let $\pts_2$ be a set of $n$ points on a curve $\gamma_2$.
We already know that there could be $\Theta(m+n)$ distinct distances when $\gamma_1$ and $\gamma_2$ are parallel or orthogonal lines.
In this more general scenario there exists a third exceptional case --- there could be $\Theta(m+n)$ distinct distances when $\gamma_1$ and $\gamma_2$ are concentric circles.
We denote by $D_{\text{curves}}(m,n)$ the minimum number of distinct distances that can occur when $\gamma_1$ and $\gamma_2$ do not contain parallel lines, orthogonal lines, and concentric circles. 
Pach and de Zeeuw \cite{PZ13} generalized \cite{SS13} to obtain $D_{\text{curves}}(m,n) = \Omega(\min\{n^{2/3}m^{2/3},m^2,n^2\})$.

\begin{problem} \label{pr:curves}
Find the asymptotic value of $D_{\text{\em curves}}(m,n)$.
\end{problem}

So far we discussed the case where both point sets are unrestricted and the case were both point sets are restricted to curves.
We can also consider the case where exactly one of the two point sets is restricted.
We denote by $D_{\text{line}}(m,n)$ the minimum number of distinct distances between a set of $m$ points on a line and a set of $n$ unrestricted points, both in $\RR^2$.
Elekes \cite{Elek95} proved that when $n\ge 4m^3$, we have $D_{\text{line}}(m,n) = O(m^{1/2}n^{1/2})$.
It is not clear whether similar constructions exist for larger values of $m$, and it is possible that when $m> (n/4)^{1/3}$ the number of distinct distances jumps to $\Omega(n/\sqrt{\log n}+m)$.

Pohoata and Sheffer \cite{PS17} derived three lower bounds for this problem: the bound $D_{\text{line}}(m,n) = \Omega(m^{1/2}n^{1/2})$ when $m=\Omega(n^{1/2}/\log^{1/3} n)$, the bound $D_{\text{line}}(m,n) = \Omega\left(n^{3/8}m^{3/4}\right)$ when $m=O(n^{1/2}/\log^{1/3} n)$ and $m=\Omega(n^{3/10})$, and the bound $D_{\text{line}}(m,n) = \Omega\left(n^{1/2}m^{1/3}\right)$ when $m=O(n^{3/10})$. 
Note that there are polynomial gaps between the lower and upper bounds in all of the above ranges. 

\begin{problem} \label{pr:OneLine}
Find the asymptotic value of $D_{\text{\em line}}(m,n)$.
\end{problem}

Let $\ell$ be the line containing $\pts_1$ in the above problem.
Bruner and Sharir \cite{BS16} studied this problem with the extra restriction that every line parallel or orthogonal to $\ell$ contains $O(1)$ points of $\pts_2$.
In this case, they proved that the number of distinct distances is
\begin{equation} \label{eq:Bruner} 
\Omega(\min\{n^{2/3}m^{2/3},m^{4/11}n^{10/11}\log^{-2/11}m,m^2,n^2\}). 
\end{equation}

\begin{problem} \label{pr:OneLineExtraRes}
Find the asymptotic value of $D_{\text{\em line}}(m,n)$ when also assuming that every line parallel or orthogonal to the line  contains $O(1)$ points of $\pts_2$.
\end{problem}

Another bipartite problem with unrestricted point sets involves $D(3,n)$. 
That is, we wish to find the minimum number of distances between $n$ points and three points. 
Recalling Elekes' bound $D_{\text{line}}(m,n) = O(m^{1/2}n^{1/2})$ from \cite{Elek95}, we obtain $D(3,n) = O(n^{1/2})$.
It is not difficult to show that $D(3,n) = \Theta(n^{1/2})$. 

The problem of $D(3,n)$ becomes more challenging when assuming that the three points are not collinear.
Elekes and Szab\'o \cite{ES12} proved that the number of distinct distances in this case is $\Omega(n^{0.502})$, showing that collinearity is necessary for obtaining a bound of $\Theta(n^{1/2})$. 
Sharir and Solymosi \cite{SSo13} improved this bound to $\Omega(n^{6/11})$. 
The current best upper bound when the points are not collinear is the trivial $D(3,n)= (n/\sqrt{\log n})$.

\begin{problem}
Find the asymptotic value of $D(3,n)$ when the three points are not collinear.
\end{problem}

\parag{Higher dimensions.}
The techniques for the above bipartite planar problems extend to some problems in higher dimensions. 
For example, Bruner and Sharir \cite{BS16} also obtained \eqref{eq:Bruner} in the case where the points of $\pts_1$ are on a line $\ell \subset \RR^d$ and every hyperplane orthogonal to $\ell$ and hypercylinder having $\ell$ as its axis contains $O(1)$ points of $\pts_2$.

As another example, consider the minimum number of distances between a set $\pts_1$ of $n$ points on a surface $S_1$ and a set $\pts_2$ of $n$ points on a surface $S_2$, both in $\RR^3$.
When $S_1$ and $S_2$ are non-parallel planes, there are two orthogonal lines $\ell_1,\ell_2$ such that $\ell_1 \subset S_1$, $\ell_2 \subset S_2$, and $\ell_1 \cap \ell_2$ is a point on $S_1 \cap S_2$. 
By placing points on $\ell_1,\ell_2$ as depicted in Figure \ref{fi:OrthLines}, we obtain $\Theta(m+n)$ distinct distances. 
The same bound can be obtained between two parallel planes, between two spheres, and between a sphere and a plane.
However, even these special cases are far from being settled, since an unrestricted set of $n$ points in $\RR^3$ can determine $O(n^{2/3})$ distinct distances.

Sharir and Solomon \cite{ShSo17} studied the case where $\pts_1$ is on a constant-degree surface in $\RR^3$ and $\pts_2$ is an unrestricted set in $\RR^3$.
In this case, they derived the bound $D(\pts_1,\pts_2) = \Omega(\min\{m^{4/7-\eps}n^{1/7},m,n\})$.

Bipartite bounds can also be obtained in various other scenarios in higher dimensions. 
It is not yet clear what the main problems and difficulties are, so we conclude this section with a deliberately vague problem.

\begin{problem} 
Derive non-trivial bounds for bipartite distinct distances problems in $\RR^d$, for $d\ge 3$.
\end{problem}

\section{Subsets with no repeated distances} \label{sec:subset}

This section surveys problems that concern point subsets that do not span a distance more than once.
Table \ref{ta:subset} lists the current best bounds for the problems that are presented in this section.

\begin{table}[h]
\begin{center}
\caption{Current best bounds for the problems of Section \ref{sec:subset}. \label{ta:subset}}
\vspace{2mm}

\begin{tabular}{|c|c|c|}
\hline
{\textbf {Variant}} & {\textbf {Lower bound}} & {\textbf {Upper bound}} \\
\hline\hline
$\s(n)$  & $\Omega(n^{1/3}/\log^{1/3} n)$ \cite{Chara12,GK15,LT95} & $O\left(\sqrt{n}/(\log n)^{1/4}\right)$ \cite{EG70} \\ \hline
\ignore{$\s_O(n)$  & $\Omega(n^{1/3}/\log^{1/3} n)$ \cite{Chara12,Tao11} & $O\left(\sqrt{n}\right)$ (trivial) \\ \hline}
$\s(\lattice)$  & $\Omega(n^{1/3}/\log^{1/3} n)$ \cite{EG70} & $O\left(\sqrt{n}/(\log n)^{1/4}\right)$ \cite{EG70} \\ \hline
$\s(\lattice_d)$  & $\Omega(n^{2/(3d)})$ \cite{LT95} & $O(n^{1/d})$ (trivial) \\ \hline
$\s_d(n)$ & $\Omega\left(n^{1/(3d-3)}(\log n)^{1/3 -2/(3d-3)}\right)$ \cite{CFGHUZ14} & $O(n^{1/d})$ (trivial) \\ \hline
$\s'(n)$  & $\Omega(n^{0.4315})$ & $O\left(\sqrt{n}/(\log n)^{1/4}\right)$ \cite{EG70} \\ \hline

\end{tabular}
\end{center}
\vspace{-\baselineskip}
\end{table}

\begin{figure}[h]

\centerline{\placefig{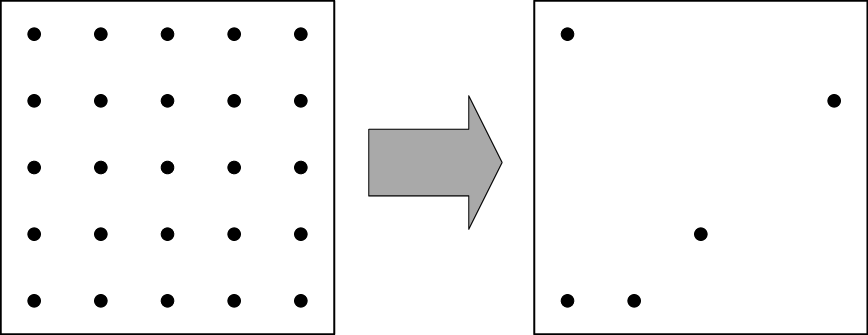}{0.5\textwidth}}
\vspace{-1mm}

\caption{\small \sf A set of 25 points and a subset of four points that span every distance at most once. No subset of five points has this property.}
\label{fi:gridSubset}
\vspace{-2mm}
\end{figure}

Given a set $\pts$ of points in $\RR^2$, let $\s(\pts)$ denote the size of the largest subset $\pts' \subset \pts$ such that every distance is spanned by the points of $\pts'$ at most once; that is, there are no points $a,b,c,d \in \pts'$ such that $|ab|=|cd|>0$ (including cases where $a=c$).
Figure \ref{fi:gridSubset} depicts a set of 25 points and a subset of four points that span every distance at most once.
Let $\s(n) = \min_{|\pts|=n}\s(\pts)$. 
In other words, $\s(n)$ is the maximum number satisfying that every set of $n$ points in $\RR^2$ contains a subset of $\s(n)$ points that do not span any distance more than once.

\begin{problem} \label{pr:subset} {\bf (Erd\H os \cite{Erd57,EG70,EP95})}
Find the asymptotic value of $\s(n)$.
\end{problem}

Let $\pts$ be a point set, such that no distance is spanned more than once by a subset $\pts'\subset \pts$. 
Then $D(\pts)\ge \binom{|\pts'|}{2}$, or equivalently $|\pts'|=O\left(\sqrt{D(\pts)}\right)$.
Let $\lattice$ be a $\sqrt{n}\times \sqrt{n}$ section of $\ZZ^2$. 
Recall from the introduction and Section \ref{sec:structure} that $D(\lattice)=\Theta(n/\sqrt{\log n})$. 
Therefore, we have $\s(n) \le \s(\lattice) = O\left(\sqrt{n}/(\log n)^{1/4}\right)$.
Lefmann and Thiele \cite{LT95} used a probabilistic argument to derive the bound $\s(n)=\Omega(n^{0.25})$. 
Dumitrescu \cite{Dum08} improved this bound to $\s(n)=\Omega(n^{0.288})$.
Charalambides \cite{Chara12} obtained the following elegant improved result by combining the probabilistic argument of Lefmann and Thiele with a result from Guth and Katz's distinct distances paper \cite{GK15}.

\begin{theorem} \label{th:subset} {\bf (Charalambides \cite{Chara12})}
$\s(n)=\Omega(n^{1/3}/\log^{1/3} n)$.
\end{theorem}
\begin{proof}
Consider a set $\pts$ of $n$ points in $\RR^2$, and define the set
\[ Q_1 = \left\{ (a,b,c,d)\in \pts^4 \ \big| \ |ab|=|cd|>0 \right\}, \]
where every quadruple of $Q_1$ consists of four distinct points. Guth and Katz \cite{GK15} proved that $|Q_1|=O(n^3\log n)$. Let $Q_2$ be the set of isosceles and equilateral triangles that are spanned by points of $\pts$.
Pach and Tardos \cite{PT02} proved that $|Q_2| = O(n^{2.137}$).

Let $\pts' \subset \pts$ be a subset that is obtained by selecting every point of $\pts$ with a probability $0 < p <1$ that will be determined below.
We have $\EE[|\pts'|] = p n$. Let $Q_1' \subset Q_1$ be the set of quadruples of $Q_1$ that contain only points of $\pts'$. Every quadruple of $Q_1$ is in $Q_1'$ with a probability of $p^4$, so $\EE[|Q_1'|] \le \alpha p^4 n^3\log{n}$, for a sufficiently large constant $\alpha$. 
Let $Q_2'$ be the set of triangles of $Q_2$ that contain only points of $\pts'$, and note that $\EE[|Q_2'|] \le \alpha p^3 n^{2.137}$ for sufficiently large $\alpha$.
Note that the points of $\pts'$ span every distance at most once if and only if $|Q_1'|=|Q_2'|=0$.
By linearity of expectation, we have
\begin{equation*}
\EE\left[|\pts'|-|Q_1'|-|Q_2'|\right]\ge p n - \alpha p^4 n^3\log{n} - \alpha p^3 n^{2.137}.
\end{equation*}
By setting $p=1/(2\alpha n^2\log n)^{1/3}$, for sufficiently large $n$ we obtain
\begin{equation*}
\EE\left[|\pts'|-|Q_1'|-|Q_2'|\right]> \frac{n^{1/3}}{3(\alpha\log n)^{1/3}}.
\end{equation*}
Therefore, there exists a subset $\pts' \subset \pts$ for which $|\pts'|-|Q_1'|-|Q_2'|\ge \frac{n^{1/3}}{3(\alpha\log n)^{1/3}}$. Let $\pts''$ be a subset of $\pts'$ that is obtained by removing from $\pts'$ an arbitrary point from every element of $Q_1'$ and $Q_2'$. The subset $\pts''$ does not span any repeated distances and contains $\Omega(n^{1/3}/\log^{1/3} n)$ points of $\pts$.
\end{proof}

When allowing a distance to repeat a small number of times, one can find larger subsets.
Pohoata and Sheffer \cite{PS17} proved that in every set of $n$ points in $\RR^2$ there exists a subset of size $\Omega(n^{22/63}\log^{-13/63}n)$ with no distance repeating more than four times.
Similarly, there exists a subset of size $\Omega(n^{12/35}\log^{-9/63}n)$ with no distance repeating more than twice. 

Erd\H os and Guy \cite{EG70} considered the following special case of Problem \ref{pr:subset}.
\begin{problem} \label{pr:Lat}{\bf (Erd\H os and Guy \cite{EG70})}
Find the asymptotic value of $\s(\lattice)$, where $\lattice$ is a $\sqrt{n}\times \sqrt{n}$ integer lattice.
\end{problem}

As mentioned above, the current best upper bound for $\s(\lattice)$ is $O\left(\sqrt{n}/(\log n)^{1/4}\right)$. Erd\H os and Guy \cite{EG70} derived the bound $\s(\lattice) = \Omega(n^{1/3-\varepsilon})$, which was later improved by Lefmann and Thiele \cite{LT95} to $\s(\lattice) = \Omega(n^{1/3})$. This bound is still marginally better than the bound implied by Theorem \ref{th:subset}.

Erd\H os and Guy \cite{EG70} also considered the higher-dimensional variant of Problem \ref{pr:Lat}. 
That is, they considered a $d$-dimensional lattice $\lattice_d$ of the form $n^{1/d}\times \cdots \times n^{1/d}$.
Erd\H os and Guy derived the bound $\s(\lattice_d) = \Omega(n^{2/(3d)-\varepsilon})$, and this was later improved by Lefmann and Thiele \cite{LT95} to $\s(\lattice_d) = \Omega(n^{2/3d})$.
It is simple to show that the points of $\lattice_d$ span $O(n^{2/d})$ distinct distances (see Section \ref{sec:HighD}), which implies $\s(\lattice_d) = O(n^{1/d})$.

\begin{problem}{\bf (Erd\H os and Guy \cite{EG70})}
Find the asymptotic value of $\s(\lattice_d)$, where $\lattice_d$ is an $n^{1/d}\times \cdots \times n^{1/d}$ integer lattice.
\end{problem}

One can also consider the higher-dimensional variant of Problem \ref{pr:subset}. 
Let $\s_d(n)$ denote the maximum number satisfying that every set of $n$ points in $\RR^d$ contains a subset of $\s_d(n)$ points that do not span any distance more than once.
Thiele \cite[Theorem 4.33]{Thiele95} proved the lower bound $\s_d(n)=\Omega(n^{1/(3d-2)})$.
This was improved by Conlon, Fox, Gasarch, Harris, Ulrich, and Zbarsky \cite{CFGHUZ14} to $\s_d(n)=\Omega\left(n^{1/(3d-3)}(\log n)^{1/3 -2/(3d-3)}\right)$.
The current best upper bound is $\s_d(n)\le \s(\lattice_d) = O(n^{1/d})$.
\begin{problem}
Find the asymptotic value of $\s_d(n)$ for $d\ge 3$.
\end{problem}

The open problems book of Brass, Moser, and Pach \cite{BMP05} offers another problem of a similar flavor. 
Let $\s'(n)$ denote the maximum number satisfying the property that every set of $n$ points in the plane contains a subset of $\s'(n)$ points that do not span any isosceles triangles.

\begin{problem} {\bf (Brass, Moser, and Pach \cite{BMP05})}
Find the asymptotic value of $\s'(n)$.
\end{problem}

For a trivial upper bound, we have $\s'(n) \le s(n) = O\left(\sqrt{n}/(\log n)^{1/4}\right)$. 
By adapting the proof of Theorem \ref{th:subset} (that is, removing $Q_1$ from the analysis), one obtains $s'(n) = \Omega(n^{0.4315})$.

\section{Distinct distances with local properties}\label{sec:kl}

For positive integers $k,\ell$, we consider planar point sets where every $k$ points determine at least $\ell$ distinct distances.
Let $\phi(n,k,\ell)$ denote the minimum number of distinct distances that are span by such a set of $n$ points. 
That is, by having a local property of every small subset of points, we wish to obtain a global property of the entire point set.
Studying $\phi(n,k,\ell)$ was originally suggested by Erd\H os (for example, see \cite{erd86}).
Table \ref{ta:secKL} lists some of the current best bounds for small values of $k$.

\begin{table}[h]
\begin{center}
\caption{The current best bounds for some of the problems that are presented in Section \ref{sec:kl}.\label{ta:secKL}}
\vspace{2mm}

\begin{tabular}{|c|c|c|}
\hline
{\textbf {Variant}} & {\textbf {Lower bound}} & {\textbf {Upper bound}} \\
\hline\hline
$\phi(n,3,2)$  & $\Omega(n/\log n)$ \cite{GK15} & $O(n/\sqrt{\log n})$ \cite{erd46} \\ \hline
$\phi(n,3,3)$  & $\Omega(n)$ & $n2^{O(\sqrt{\log n})}$ \cite{erd86} \\ \hline
$\phi(n,4,3)$  & $\Omega(n/\log n)$ \cite{GK15} & $O(n/\sqrt{\log n})$ \\ \hline
$\phi(n,4,4)$  & $\Omega(n/\log n)$ \cite{GK15} & $n2^{O(\sqrt{\log n})}$ \cite{Dum08} \\ \hline
$\phi(n,4,5)$  & $\Omega(n)$ (trivial) & $O(n^2)$ (trivial) \\ \hline
$\phi(n,5,9)$  & $\Omega(n)$ (trivial) & $O(n^2)$ (trivial) \\ \hline
\end{tabular}
\end{center}
\vspace{-\baselineskip}
\end{table}

\parag{The case of $\mathbf{\emph{k=3}}$.} The value of $\phi(n,3,2)$ is the minimum number of distinct distances determined by a set of $n$ points that do not span any equilateral triangles. As discussed in Section \ref{sec:structure}, Erd\H os \cite{erd46} noticed that a $\sqrt{n}\times\sqrt{n}$ integer lattice determines $\Theta(n/\sqrt{\log n})$ distinct distances. It is known that the points of the integer lattice do not determine any equilateral triangles, and thus $\phi(n,3,2)= O(n/\sqrt{\log n})$.
Guth and Katz's bound implies $\phi(n,3,2) \ge D(n)= \Omega(n/\log n)$. 
Thus, the current best bounds for $\phi(n,3,2)$ are identical to the ones for $D(n)$.

\begin{problem} \label{pr:n32}
Find the asymptotic value of $\phi(n,3,2)$.
\end{problem}

The value of $\phi(n,3,3)$ is the minimum number of distinct distances determined by a set of $n$ points that do not span any isosceles triangles. 
Here and in the following cases we also consider degenerate polygons, such as a degenerate isosceles triangle whose three vertices are collinear. 
Since no isosceles triangles are allowed, a point $p\in \pts$ determines $n-1$ distinct distances with the points of $\pts\setminus\{p\}$, so $\phi(n,3,3) = \Omega(n)$.
Erd\H os \cite{erd86} observed the following upper bound for $\phi(n,3,3)$.
Behrend \cite{Behr46} proved that there exists a set $A$ of positive integers $a_1<a_2<\cdots<a_n$, such that no three elements of $A$ determine an arithmetic progression and $a_n < n2^{O(\sqrt{\log n})}$.
Note that the point set $\pts_1 = \{(a_1,0), (a_2,0),\ldots, (a_n,0)\} \subset \RR^2$ does not span any isosceles triangles. 
Since $\pts_1 \subset \pts_2 = \{(1,0),(2,0),\ldots,(a_n,0) \}$ and $D(\pts_2)< n2^{O(\sqrt{\log n})}$, we have $\phi(n,3,3) < n2^{O(\sqrt{\log n})}$.
Erd\H os conjectured \cite{erd86} that $\phi(n,3,3) = \omega(n)$.

\begin{problem} \label{pr:n33}
Find the asymptotic value of $\phi(n,3,3)$.
\end{problem}

\parag{The case of $\mathbf{\emph{k=4}}$.} The value of $\phi(n,4,3)$ is the minimum number of distinct distances determined by a set of $n$ points that do not span any squares. 
The $\sqrt{n}\times\sqrt{n}$ triangular lattice determines $\Theta(n/\sqrt{\log n})$ distinct distances (for example, see \cite{Sblog14}). Since the triangular lattice does not contain any squares, we have $\phi(n,4,3) = O(n/\sqrt{\log n})$. 
The current best lower bound is $\phi(n,4,3) \ge D(n) = \Omega(n/\log n)$.

\begin{problem} \label{pr:n43}
Find the asymptotic value of $\phi(n,4,3)$.
\end{problem}

\begin{figure}[h]

\centerline{\placefig{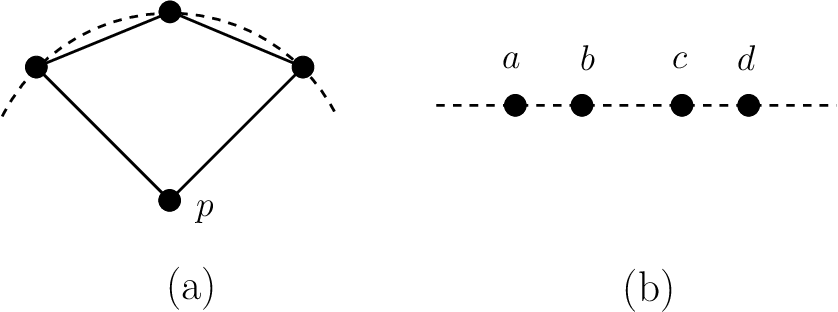}{0.6\textwidth}}
\vspace{-1mm}

\caption{\small \sf (a) The point $p$ is equidistant from the other three vertices of the deltoid. (b) The only possible cases of segments having the same length are $|ab|=|cd|$ and $|ac|=|bd|$.}
\label{fi:deltoid}
\vspace{-2mm}
\end{figure}

The value of $\phi(n,4,4)$ is the minimum number of distinct distances determined by a set of $n$ points that do not span any rhombuses, rectangles, or deltoids with one vertex that is equidistant to the three other three (see Figure \ref{fi:deltoid}(a)). Dumitrescu \cite{Dum08} observed that $\phi(n,4,4)< n2^{O(\sqrt{\log n})}$ by using the same point set $\pts_1$ from the analysis of $\phi(n,3,3)$. Indeed, consider a subset of four points of $\pts_1$, as depicted in Figure \ref{fi:deltoid}(b), and notice that the only pairs of segments that are allowed to have the same length (without resulting in an arithmetic progression) are $|ab|=|cd|$ and $|ac|=|bd|$. Thus, every quadruple of points determines at least four distinct distances.
No lower bound is known beyond $\phi(n,4,4)=\Omega(n/\log n)$.

\begin{problem} \label{pr:n44}
Find the asymptotic value of $\phi(n,4,4)$.
\end{problem}

Not much is known about the case of $\phi(n,4,5)$, even though it is considered to be one of the main variants of the problem. 
While Erd\H os \cite{erd86} asked whether $\phi(n,4,5)=\Theta(n^2)$, the current best lower bound is only $\phi(n,4,5)= \Omega(n)$.
Indeed, in a set of $n$ point set $\pts$ with this property, any circle whose center is a point of $\pts$ can be incident to at most two points of $\pts$.

\begin{problem} \label{pr:n45}
Find the asymptotic value of $\phi(n,4,5)$.
\end{problem}

We also note the trivial bound $\phi(n,4,6) = \binom{n}{2}$, since every distance can occur at most once in this case.

\parag{The case of $\mathbf{\emph{k=5}}$.} 
We do not go over the various cases of $k=5$, but only mention one that is considered interesting.
Erd\H os \cite{erd86} asked whether $\phi(n,5,9)=\Omega(n^2)$. 
Unfortunately, nothing is known in this case beyond the trivial $\phi(n,5,9)=\Omega(n)$.

\begin{problem} \label{pr:n59}
Find the asymptotic value of $\phi(n,5,9)$.
\end{problem}

\parag{Larger values of $k$.} 
We now discuss bounds for $\phi(n,k,\ell)$ that hold for arbitrarily large $k$.
A first simple observation is that for any $k\ge 4$, we have
\[ \phi\left(n,k,\binom{k}{2}-\lfloor k/2 \rfloor +2 \right) = \Omega(n^2). \]
Indeed, with this restriction every distance can occur at most $\lfloor k/2 \rfloor -1$ times.
Since every distance repeats a constant number of times, there must be a quadratic number of distances.

Similarly to the cases of $\phi(n,3,3)$ and $\phi(n,4,4)$, we can use the set $\pts_1$ to obtain the upper bound $n2^{O(\sqrt{\log n})}$ for various $\phi(n,k,\ell)$. For example, it is not difficult to show that
\[\phi\left(n,k,2\lfloor k/2 \rfloor\right)=n2^{O(\sqrt{\log n})}. \]
It seems likely that a more careful analysis would yield the same bound for larger values of $\ell$.

For any $\eps>0$, Fox, Pach, and Suk \cite{FPS18} derived the bound
\begin{equation} \label{eq:FPSlocal}
\phi\left(n,k,\binom{k}{2} - k + 6\right)=\Omega\left(n^{8/7-\eps}\right). 
\end{equation}
For any integers $k > m \ge 2$, Pohoata and Sheffer \cite{PS17} derived the bound
\[ \phi\left(n,k,\binom{k}{2} -  m\cdot \left\lfloor\frac{k}{m+1}\right\rfloor  + m+1 \right) = \Omega\left(n^{1+\frac{1}{m}}\right). \]
This bound is stronger than \eqref{eq:FPSlocal} when $m\le 7$.
For example, when $m=2$ it leads to 
\[ \phi\left(n,k,\binom{k}{2} -  2\cdot \left\lfloor\frac{k}{3}\right\rfloor  + 3 \right) = \Omega\left(n^{3/2}\right). \]

This is the only non-trivial distinct distances problem for which we can derive a bound asymptotically larger than $n^{4/3}$.
However, even the stronger $\Omega\left(n^{3/2}\right)$ is likely to be far from optimal.
It is plausible that the correct values in these cases is close to $n^2$.

\begin{problem} \label{pr:nkl}
Find stronger general bounds for $\phi(n,k,\ell)$.
\end{problem}

\section{Problems related to Additive Combinatorics} \label{sec:additive}

Connections between the fields of Discrete Geometry and Additive Combinatorics are constantly being discovered. 
The purpose of this section is to mention some of these connections that involve distinct distances problems.

Given a finite set $A\subset \RR^d$, the \emph{difference set} of $A$ is defined as 
\[ A-A = \left\{a-a' :\ a,a'\in A \right\}. \]
The additive energy of $A$ is defined as
\[ E(A) = |\{(a,b,c,d)\in A^4 :\, a+b=c+d\}|. \]
A standard reference for further reading about these concepts is the book of Tao and Vu \cite{TV06}.

Nets Katz stated the following problem, suggesting that it might be an approach for characterizing planar point sets that span few distinct distances (Problem \ref{pr:Structure}).

\begin{problem} \label{pr:KatzConj}
Prove or disprove: If is a set $\pts$ of $n$ points in $\RR^2$ spans $O(n/\sqrt{\log n})$ distinct distances, then $E(\pts)$ is large.
\end{problem}

For a finite set $A \subset \RR$, Hanson \cite{Hanson17} studied the case where the cartesian product $A \times A \subset \RR^2$ spans few distinct distances.
Hanson's results were pushed further by Roche--Newton \cite{RN16} and by Pohoata \cite{Pohoata16}.
In particular, Pohoata proved that 
\[ |A-A| = O\left(D(A\times A)^{6/7}\log^{1/7}|A|\right), \] 
and this is the current best bound for the problem.
That is, if a cartesian product $A\times A$ spans a small number of distinct distances, then the difference set $A-A$ cannot be too large.

\begin{problem} \label{pr:Hanson}
For a finite set $A \subset \RR$, find the asymptotic value of $|A-A|$ when $A\times A$ spans few distinct distances.
\end{problem}

By recalling the constructions mentioned in Section \ref{sec:structure}, we note that there exist sets that span few distinct distances and are of the form $A\times B$ rather than $A\times A$.
From this perspective, Problem \ref{pr:Hanson} is more interesting when considered for $A\times B$. 
No non-trivial bounds are known in this bipartite case.

Pham, Vinh, and de Zeeuw \cite{PVdZ17} derived another related result: For every finite $A \subset \RR$ and integer $d\ge 2$, the $d$-dimensional cartesian product $A\times \cdots \times A \subset \RR^d$ satisfies
\[ D(A\times \cdots \times A) =\Omega\left(|A|^2\log^{-1/2^{d-2}}|A|\right). \]

\section{Additional problems} \label{sec:additional}
In this final section we discuss problems that did not fit into any of the previous topics.

A more general variant of Problem \ref{pr:main} asks to show that in every planar point set there exists a point that spans many distinct distances (for example, see Erd\H os \cite{erd87}).
Let $\hd(n)$ be the minimum value such that for every set $\pts$ of $n$ points in $\RR^2$ there exists a point $p \in \pts$ satisfying  $D(\{p\},\pts\setminus \{p\})\ge \hd(n)$.
A simple upper bound is $\hd(n) \le D(n) = O(n/\sqrt{\log n})$.
However, Guth and Katz's bound does not immediately imply a matching lower bound for $\hd(n)$.
The current best lower bound, obtained by Katz and Tardos \cite{KT04}, is $\hd(n) = \Omega(n^{(48-14e)/(55-16e)})\approx \Omega(n^{0.864})$.

\begin{problem} \label{pr:mainPnt}
Find the asymptotic value of $\hd(n)$.
\end{problem}

Given a planar point set $\pts$ and a point $p \in\pts$, we denote by $\hd_p(\pts)$ the number of distinct distances between $p$ and the other points of $\pts$. We set $\hd_\Sigma(\pts)=\sum_{p\in\pts}\hd_p(\pts)$ and $\hd_\Sigma(n) = \min_{|\pts|=n}\hd_\Sigma(\pts)$. 
It is not hard to verify that when $\pts$ is a $\sqrt{n}\times\sqrt{n}$ section of $\ZZ^2$, we have $\hd_\Sigma(\pts) = \Theta\left(n^2/\sqrt{\log n}\right)$.
This implies that $\hd_\Sigma(n) = O\left(n^2/\sqrt{\log n}\right)$, which Erd\H os conjectured \cite{Erd75,erd87} to be tight.

\begin{problem}
Find the asymptotic value of $\hd_\Sigma(n)$.
\end{problem}
The current best lower bound $\hd_\Sigma(n) = \Omega(n^{1.864})$ is immediately implied by Katz and Tardos' \cite{KT04} bound $\hd(n) = \Omega(n^{(48-14e)/(55-16e)})\approx \Omega(n^{0.864})$ (by repeatedly removing a point $p$ that maximizes $\hd_p(\pts)$).

Sharir, Smorodinsky, Valculescu, and de Zeeuw \cite{SSmo14} studied a distinct distances problem between points and lines.
Denote by $L(m,n)$ the minimum number of distinct distances between a set of $m$ points and a set of $n$ \emph{lines}, both in $\RR^2$ (where the distance between a point and a line is defined in the usual way).
By placing $m$ points on a line $\ell$ and then taking $n$ lines that are parallel to $\ell$, we obtain  $L(m,n) = O(n)$. 
The bound $L(m,n) = \Omega( m^{1/5} n^{3/5} )$ was derived in \cite{SSmo14} when $\sqrt{m}< n <m^2$.

\begin{problem}
Find the asymptotic value of $L(m,n)$.
\end{problem}

Finally, we consider a generalization from distinct distances to \emph{distinct vectors}. 
We say that points $p,q\in \RR^2$ span the vectors $p-q$ and $q-p$.
Given a planar point set $\pts$, we denote by $v(\pts)$ the number of distinct vectors that are spanned by pairs of points of $\pts$, and set $v(n) = \min_{|\pts|=n}v(\pts)$. 
It is also not difficult to show that $v(n) = \Theta(n)$. 
To obtain a lower bound, we note that each point determines at least $n-1$ distinct vectors.
To obtain an upper bound, we may consider a $\sqrt{n}\times\sqrt{n}$ section of $\ZZ^2$, or evenly spaced points on a line.

Erd\H os, F\"uredi, Pach, and Ruzsa \cite{EFPR93} studied the case of distinct vectors for point sets in general position (that is, no three points on a line and no four points on a circle). We set $v_\text{gen}(n) = \min_{|\pts|=n}v(\pts)$, where the sum is taken over every set of $n$ points in general position. In \cite{EFPR93}, it is proven that $v_\text{gen}(n) > cn$ for \emph{every} constant $c$; i.e., $v_\text{gen}(n) =\omega(n)$. The current best upper bound  $v_\text{gen}(n) = n2^{O(\log n)}$ is immediately implied by the current best upper bound for Problem \ref{pr:gen}.

\begin{problem}
Find the asymptotic value of $v_\text{gen}(n)$.
\end{problem}



\begin{thebibliography}{99}
%
\bibitem{Alt63}
E.\ Altman,
On a problem of P.\ Erd\H os,
\emph{Amer.\ Math.\ Monthly} {\bf 70} (1963), 148--157.
%
\bibitem{Alt72}
E.\ Altman,
Some theorems on convex polygons,
\emph{Canad.\ Math.\ Bull.} {\bf 15} (1972), 329--340.
%
\bibitem{BES17}
S.\ Bardwell-Evans and A.\ Sheffer,
A Reduction for the Distinct Distances Problem in $\RR^d$,
arXiv:1705.10963.
%
\bibitem{Behr46}
F.\ A.\ Behrend,
On sets of integers which contain no three terms in arithmetic progression,
\emph{Proc.\ Nat.\ Acad.\ Sci.} {\bf 32} (1946), 331--332.
%
\bibitem{BR95}
B.\ C.\ Berndt and R.\ A.\ Rankin,
\emph{Ramanujan: Letters and Commentary},
Amer. Math. Soc., Providence, RI, 1995.
%
\bibitem{BMP05}
P.\ Brass, W.\ Moser, and J.\ Pach,
\emph{Research Problems in Discrete Geometry},
Springer-Verlag, New York, 2005.
%
\bibitem{BMO11}
D.\ Brink, P.\ Moree, and R.\ Osburn,
Principal forms $X^2 + nY^2$ representing many integers,
\emph{Abh.\ Math.\ Sem.\ Univ.\ Hambg.}, {\bf 81} (2011), 129--139.
%
\bibitem{BS16}
A.\ Bruner and M.\ Sharir,
Distinct distances between a collinear set and an arbitrary set of points,
\emph{Discrete Math.}, {\bf 341} (2018), 261--265.
%
\bibitem{Chara12}
M.\ Charalambides,
A note on distinct distance subsets,
\emph{Journal of Geometry}, {\bf 104} (2013), 439--442.
%
\bibitem{Chara14}
M.\ Charalambides, 
Distinct distances on curves via rigidity, 
\emph{Discrete Comput. Geom.}, {\bf 51} (2014), 666--701.
%
\bibitem{CFGHUZ14}
D.\ Conlon, J.\ Fox, W.\ Gasarch, D.\ Harris, D.\ Ulrich, and S.\ Zbarsky,
Distinct volume subsets,
\emph{SIAM Journal on Discrete Mathematics}, {\bf 29} (2015), 472--480.
%
\bibitem{Dum06}
A.\ Dumitrescu,
On distinct distances from a vertex of a convex polygon,
\emph{Discrete Comput. Geom.} {\bf 36} (2006), 503--509.
%
\bibitem{Dum08}
A.\ Dumitrescu, On distinct distances among points in general position and other related problems,
\emph{Periodica Mathematica Hungarica} {\bf 57} (2008), 165--176.
%
\bibitem{Elek95}
G.\ Elekes,
Circle grids and bipartite graphs of distances,
\emph{Combinatorica} {\bf 15} (1995), 167--174.
%
\bibitem{Elek99}
G.\ Elekes,
A note on the number of distinct distances,
\emph{Period.\ Math.\ Hung.}, {\bf 38} (1999), 173--177.
%
\bibitem{ER00}
G.\ Elekes and L.\ R\'onyai, A combinatorial problem on polynomials and
rational functions,
\emph{J.\ Combin.\ Theory Ser.\ A}, {\bf 89} (2000), 1--20.
%
\bibitem{ES12}
G.\ Elekes and E.\ Szab\'o,
How to ﬁnd groups? (and how to use them in Erd\H os geometry?),
\emph{Combinatorica} {\bf 32} (2012), 537--571.
%
\bibitem{erd46}
P.\ Erd\H os,
On sets of distances of $n$ points,
\emph{Amer. Math. Monthly} {\bf 53} (1946), 248--250.
%
\bibitem{Erd57}
P.\ Erd\H os, Nehany geometriai probl\'em\'ar\'ol (in Hungarian), \emph{Matematikai Lapok}, {\bf 8} (1957), 86--92.
%
\bibitem{Erd75}
P.\ Erd\H os,
On some problems of elementary and combinatorial geometry,
\emph{Ann.\ Mat.\ Pura Appl.} {\bf 103} (1975), 99--108.
%
\bibitem{erd86}
P.\ Erd\H os,
On some metric and combinatorial geometric problems,
\emph{Discrete Math.} {\bf 60} (1986), 147--153.
%
\bibitem{erd87}
P.\ Erd\H os,
Some combinatorial and metric problems in geometry,
\emph{Intuitive geometry} (K.\ B\"or\"oczky and G.\ Fejes T\`oth, eds.), North-Holland, Amsterdam-New York, 1987, 167--177.
%
\bibitem{erd96}
P.\ Erd\H os,
On some of my favourite theorems,
\emph{Combinatorics, Paul Erd\H os is Eighty,}
Vol. 2 (D. Mikl\'os et al., eds.), Bolyai Society Mathematical Studies 2, Budapest, 1996,
97--132.
%
\bibitem{EF96}
P.\ Erd\H os and P.\ Fishburn,
Maximum planar sets that determine $k$ distances,
\emph{Discrete Math.} {\bf 160} (1996), 115--125.
%
\bibitem{EFPR93}
P.\ Erd\H os, Z.\ F\"uredi, J.\ Pach, and I.\ Z.\ Ruzsa,
The grid revisited,
\emph{Discrete Math.} {\bf 111} (1993), 189--196.
%
\bibitem{EG70}
P.\ Erd\H os, R.\ K.\ Guy,
Distinct distances between lattice points,
\emph{Elemente Math.} {\bf 25} (1970) 121--123.
%
\bibitem{EHP89}
P.\ Erd\H os, D.\ Hickerson, and J.\ Pach,
A problem of Leo Moser about repeated distances on the sphere,
\emph{Amer.\ Math.\ Monthly} {\bf 96} (1989), 569--575.
%
\bibitem{EP95}
P.\ Erd\H os and G.\ Purdy,
Extremal problems in combinatorial geometry,
\emph{Handbook of Combinatorics},
Vol. I (R.\ L.\ Graham, M.\ Gr\H otschel, and L. Lov\'asz, editors), Elsevier, Amsterdam, 1995, pp. 809--
874.
%
\bibitem{FPS18}
J.\ Fox, J.\ Pach, and A.\ Suk. 
More distinct distances under local conditions,
\emph{Combinatorica} {\bf 38} (2018), 501--509.
%
\bibitem{GIS11}	
J.\ Garibaldi, A.\ Iosevich, and S.\ Senger,
\emph{The Erd{\H o}s Distance Problem},
Student Math. Library, Vol. 56,
Amer. Math. Soc. Press, Providence, RI, 2011.
%
\bibitem{GK15}	
L.\ Guth and N.\ H.\ Katz,
On the Erd{\H o}s distinct distances problem in the plane,
{\em Annals of Mathematics}, {\bf 181} (2015), 155--190.
%
\bibitem{Hanson17}
B.\ Hanson,
The additive structure of cartesian products spanning few distinct distances, 
\emph{Combinatorica}, to appear.
%
\bibitem{KT04}
N.\ H.\ Katz and G.\ Tardos,
A new entropy inequality for the Erd\H os distance problem,
\emph{Towards a Theory of Geometric Graphs (J. Pach, ed.)}, Contemporary Mathematics {\bf 342}, AMS, Providence, RI, 2004, 119--126.
%
\bibitem{Land08}
E.\ Landau,
\"Uber die Einteilung der positiven ganzen Zahlen in vier Klassen nach der Mindeszahl der zu ihrer additiven Zusammensetzung erforderlichen Quadrate,
\emph{Arch.\ Math.\ Phys.} {\bf 13} (1908), 305--312.
%
\bibitem{LT95}
H.\ Lefmann and T.\ Thiele,
Point sets with distinct distances,
\emph{Combinatorica} {\bf 15} (1995) 379--408.
%
\bibitem{LSdZ15}
B.\ Lund, A.\ Sheffer, and F.\ de Zeeuw,
Bisector energy and few distinct distances,
\emph{Discrete Comput.\ Geom.}, {\bf 56} (2016), 337--356.
%
\bibitem{NPPZ11}
G.\ Nivasch, J.\ Pach, R.\ Pinchasi, and S.\ Zerbib,
The number of distinct distances from a vertex of a convex polygon,
\emph{J.\ Computational Geometry}, {\bf 4} (2013), 1--12.
%
\bibitem{PZ13}
J.\ Pach and F.\ de Zeeuw
Distinct Distances on Algebraic Curves in the Plane,
\emph{Proc.\ 30th ACM Symp.\ on Computational Geometry} (2014), 549--557.
%
\bibitem{PT02}
J.\ Pach and G.\ Tardos,
Isosceles triangles determined by a planar point set,
\emph{Graphs and Combinatorics} {\bf 18} (2002), 769--779.
%
\bibitem{PVdZ17}
T.\ Pham, L.\ A.\ Vinh, and F.\ de Zeeuw,
Three-variable expanding polynomials and higher-dimensional distinct distances, 
\emph{Combinatorica}, to appear.
%
\bibitem{Pohoata16}
C.\ Pohoata, 
An Update on Cartesian Products with Few Distinct Distances, 
arXiv:1612.06153.
%
\bibitem{PS17}
C.\ Pohoata and A.\ Sheffer,
Higher Distance Energies and Expanders with Structure, 
arXiv:1709.06696.
%
\bibitem{Raz16}
O.\ E.\ Raz, 
A note on distinct distances, 
arXiv:1603.00740.
%
\bibitem{RRNS15}
O.\ E.\ Raz, O.\ Roche-Newton, and M.\ Sharir,
Sets with few distinct distances do not have heavy lines,
\emph{Discrete Math.} {\bf 338} (2015), 1484--1492.
%
\bibitem{RSS14}
O.\ E.\ Raz, M.\ Sharir and J.\ Solymosi,
Polynomials vanishing on grids: The Elekes-R\'onyai problem revisited,
\emph{Amer. J. Math.}, {\bf 138} (2016), 1029--1065.
%
\bibitem{Reid70}
C.\ Reid,
\emph{Hilbert: with an appreciation of Hilbert's mathematical work by Hermann Weyl},
Springer, 1970.
%
\bibitem{RN16}
O.\ Roche--Newton, 
On sets with few distinct distances, 
arXiv:1608.02775.
%
\bibitem{SS13}
M.\ Sharir, A.\ Sheffer, and J.\ Solymosi,
Distinct distances on two lines,
\emph{J.\ Combinat.\ Theory A} {\bf 120} (2013), 1732--1736.
%
\bibitem{SSmo14}
M.\ Sharir, S.\ Smorodinsky, C.\ Valculescu, and F.\ de Zeeuw,
On Distinct Distances Between Points and Lines,
\emph{Comput.\ Geom.\ Theory Appls.} {\bf 69} (2018), 2--15.
%
\bibitem{ShSo17}
M.\ Sharir and N.\ Solomon,
Incidences between points and surfaces and points and curves, and distinct and repeated distances in three dimensions,
\emph{Proc.\ 28th ACM-SIAM Symp.\ on Discrete Algorithms} (2017), 2456--2475. 
%
\bibitem{SSo13}
M.\ Sharir and J.\ Solymosi,
Distinct distances from three points,
\emph{Combinat. Probab. Comput.}, {\bf 25} (2016), 623--632.
%
\bibitem{Sblog14}
A.\ Sheffer,
Point Sets with Few Distinct Distances, blog post,
\url{https://adamsheffer.wordpress.com/2014/07/16/point-sets-with-few-distinct-distances/}
%
\bibitem{Sblog14b}
A.\ Sheffer,
Few Distinct Distances Implies Many Points on a Line, blog post,
\url{https://adamsheffer.wordpress.com/2014/10/07/few-distinct-distances-implies-many-points-on-a-line/}
%
\bibitem{SZZ13}
A.\ Sheffer, J.\ Zahl, and F.\ de Zeeuw,
Few distinct distances implies no heavy lines or circles,
\emph{Combinatorica} {\bf 36} (2016), 349--364.
%
\bibitem{SV08}
J.\ Solymosi and V.\ H.\ Vu,
Near optimal bounds for the Erd\H os distinct distances problem in high dimensions.
\emph{Combinatorica} {\bf 28} (2008), 113--125.
%
\bibitem{Tao11}
T.\ Tao,
Lines in the Euclidean group SE(2),
blog post, \url{http://terrytao.wordpress.com/2011/03/05/lines-in-the-euclidean-group-se2/}
%
\bibitem{TV06}
T.\ Tao and V.\ H.\ Vu,
\emph{Additive combinatorics}, Cambridge University Press, 2006.
%
\bibitem{Thiele95}
T.\ Thiele,
\emph{Geometric selection problems and hypergraphs}, PhD thesis, Instut fur
Mathematik II Freir Universitat Berlin, 1995.
%
\end{thebibliography}
\end{document}